\documentclass[a4paper, 12pt]{article}
\usepackage[a4paper, left=2cm, right=1cm]{geometry}
\usepackage{multirow}
\usepackage{array}
\usepackage{color}
\usepackage{amsmath,amssymb,latexsym}
\usepackage[mathscr]{eucal}
\usepackage{arydshln}
\newcommand{\beq}{\begin{equation}}
\newcommand{\eeq}{\end{equation}}
\newcommand{\bea}{\begin{eqnarray}}
\newcommand{\eea}{\end{eqnarray}}
\newcommand{\bean}{\begin{eqnarray*}}
\newcommand{\eean}{\end{eqnarray*}}

\newenvironment{proof}%
{\par\noindent\emph{Proof:\ }}%
{\ \hfill ~\rule{2mm}{2mm}\par\bigskip}
{\par\noindent\textbf{Remark:\ }}%
{\ \hfill \par\bigskip}

 \def\t{{\mathbf t}}
 
  \def\D{{\mathbf D}}

 \def\hx{h}

\makeatletter
\newcolumntype{"}{@{\hskip\tabcolsep\vrule width 2pt\hskip\tabcolsep}}
\def\hlinewd#1{%
  \noalign{\ifnum0=`}\fi\hrule \@height #1 \futurelet
   \reserved@a\@xhline}
\makeatother
\def\t{{\mathbf t}}

 \def\x{{\mathbf x}}

 \def\F{{\mathbf F}}

\date{July 2011}
\title{{
Centers of quasi-homogeneous polynomial planar systems}}

\author{
A. Algaba, \ N. Fuentes,\footnote{Corresponding author. E-mail addres: natalia.fuentes@dmat.uhu.es} \  C. Garc\'{\i}a
\\
Dept. Integrated Sciences. Faculty of Experimental Sciences\\ University of Huelva,
Spain
}
\begin{document}
\bibliographystyle{plain}

\maketitle

\abstract{{
In this paper we determine the centers of quasi-homogeneous polynomial planar vector
fields of degree $0$, $1$, $2$, $3$ and $4$. In addition, in every case we make a study of the reversibility
and the analytical integrability of each one of the above centers.Wefind polynomial centers
which are neither orbitally reversible nor analytically integrable, this is a new scenario in
respect to the one of non-degenerate and nilpotent centers.
}

\section{Introduction} 

In the study of the planar systems, one of the classic problems in
the qualitative theory of the analytical systems, is the study of
the phase portrait in a neighborhood of a singular point and, in
particular, to characterize when the singular point is a center or
a focus. A singular point is called center-focus type, also called
monodromic, when the orbits rotate around the singular point, i.e.
it is a center or a focus in the analytic case. Once established the monodromy, the
center problem determines, by studying the Ponincar\'{e} map,
when, all the orbits in a neighborhood of the singular point are
closed. This chapter deals with the classification of the centers
for a class of polynomial differential system with null linear
part (degenerate centers). The classification of the centers of
polynomial differential system with linear part $(-y,x)^T$ started
with the works of Dulac [1], Bautin [2], Kapteyn
[3,4] and Zoladek [5] for quadratic system, 
and continued with the works of Sibirskii [6] and Zoladek [7] 
for symmetric cubic system. A lot of work has been performed about 
non-degenerate centers (see [8]), but
only partial results have been reached and we are very far to
obtain a complete classification of all non-degenerate centers for
the polynomial differential systems of degree greater than or
equal to three. The nilpotent centers (i.e.,with linear part
$(y,0)^T$) are characterized theoretically (see [9-11]) but
only a few families of nilpotents centers are known (see 12-14). However, the case of degenerate centers (i.e., with
null linear part) is not characterized theoretically (see [15]).

Another problem related to the previous one is the problem of the
reversibility, i.e., when a vector field is invariant to an
involution in the state variables and the change of sign in the
time variable. (see [9, 10, 16, 17])

Finally, the third problem we study in this chapter is the
analytical integrability of a quasi-homogeneous vector field,
more specifically, to determine when a planar vector field has an
analytical first integral, i.e., a function that remains constant
along the trajectories of the system. (see[18-23])

We consider system
\begin{equation}\label{SisNpert}
\dot{\x}=\F_r(\x).
\end{equation}
with $\F_r= (P, Q)$, with $P$ and $Q$ coprime and $\F_r$ a
quasi-homogeneous vector field. In this case, it is known that the
origin is the unique (real and finite) singularity of $\F_r$.
Therefore, if $\F_r$ has a center at the origin its period annulus
is $\mathbb{R}^2\setminus\{0\}$, i.e., the center is global.

Our two motivations for the study of the quasi-homogeneous centers
are  by one hand, Theorem 4 which assures us that: \emph{a
necessary condition so that the origin of a perturbation of a
quasi-homogeneous vector field be a center is that the origin of
the quasi-homogeneous vector field be a center}. By other hand, a
recent work Llibre and Pessoa [22], which shows a
classification of these centers up to fourth degree (third degree
for us). Notice that our definition of degree of a
quasi-homogeneous vector field disagrees with the one given in
Llibre and Pessoa [22] in one unit.

In this paper, by using others techniques, we extend this study up
to fifth degree (fourth degree for us), and besides, we
characterize the integrability and the reversibility of each one
of the centers found. It is know that all non-degenerate centers
are reversible and analytically integrable (see [24,25]). Nilpotent centers are orbitally reversible (see
[9]) but there exist nilpotent centers
that are not analytically integrable (see[10]). In
this chapter is shown a new situation, we find polynomial centers which are neither orbitally reversible nor analytically integrable
(see Theorems 13 and 18). 

This paper is structured as follows. In the next section we characterize the monodromy, center and reversibility problem
for quasi-homogeneous vector fields. In Section 3 we calculate, by applying the results obtained in the previous section,
all the centers for quasi-homogeneous vector fields up to four degree, extending the results reached in [22]. We also study
which of them are reversible and analytically integrable (integrable), detecting cases of centers which are neither reversible
nor integrable.

\section{Monodromy and reversibility of quasi-homogeneous vector field}

We start this section giving some definitions and previous results about the quasi-homogeneous vector field.

\begin{itemize}
    \item Let $\mathbf{t} = (t_1, t_2)$ be non-null, with $t_1$ and $t_2$ being non-negative integer coprime numbers and $t_1 \leq t_2$. A function $p : \mathbb{R}^2 \to \mathbb{R}$ is \textit{quasi-homogeneous} of type $\mathbf{t}$ and degree $k$ if $p(\varepsilon^{t_1}x, \varepsilon^{t_2}y) = \varepsilon^k p(x, y)$. The vector space of quasi-homogeneous polynomials of type $\mathbf{t}$ and degree $k$ will be denoted by $P^t_k$.
    
    \item A polynomial vector field $\mathbf{F} = (P, Q)^T$ is quasi-homogeneous of type $\mathbf{t} = (t_1, t_2)$ and degree $r$ (denoted as $\text{deg}(\mathbf{F}) = r$) if $P \in P^t_{r + t_1}$ and $Q \in P^t_{r + t_2}$. We will denote by $Q^t_r$ the vector space of quasi-homogeneous polynomial vector fields of type $\mathbf{t}$ and degree $r$. (Notice that $\text{deg}(\mathbf{F}) = r - 1$ in the definition given in [14]).
    
    \item Let $\mathbf{t} = (t_1, t_2)$, we denote $|\mathbf{t}| = t_1 + t_2$.
    
    \item We denote by $\mathbf{X}_h := \left( -\frac{\partial h}{\partial y}, \frac{\partial h}{\partial x} \right)^T$ the Hamiltonian vector field with Hamilton function $h$, and by $\mathbf{D}_0 = (t_1 x, t_2 y)^T$ the Euler vector field associated to the type $\mathbf{t}$.
    
    \item Given a vector field $\mathbf{F} = (P, Q)^T$, we define the divergence of $\mathbf{F}$ as:
    \[
    \text{div}(\mathbf{F}) := \frac{\partial P}{\partial x} + \frac{\partial Q}{\partial y}.
    \]
    
    \item We define the wedge product of two vector fields, $\mathbf{F} \wedge \mathbf{G} := P Q' - Q P'$, where $\mathbf{F} = (P, Q)^T$ and $\mathbf{G} = (P', Q')^T$.
    
\end{itemize}

Given a type $\mathbf{t}$, any quasi-homogeneous vector field can be uniquely decomposed as the sum of two quasi-homogeneous vector fields: one with zero divergence (the \textit{conservative part}) and the other with divergence equal to the original vector field (the \textit{dissipative part}). The proof of this result can be found in [26].

\medskip

\textbf{Lemma 1.} Let $\mathbf{F}_r \in Q^t_r$. Then there exist unique polynomials $\mu \in P^t_r$ and $h \in P^t_{r + |\mathbf{t}|}$ such that:
\[
\mathbf{F}_r = \mathbf{X}_h + \mu \mathbf{D}_0,
\]
where:
\[
\mu = \frac{1}{r + |\mathbf{t}|} \, \text{div}(\mathbf{F}_r) \quad \text{and} \quad h = \frac{1}{r + |\mathbf{t}|} [\mathbf{D}_0 \wedge \mathbf{F}_r].
\]

The above expression is called the \textit{conservative–dissipative decomposition} of $\mathbf{F}_r \in Q^t_r$, where $\mathbf{X}_h$ and $\mu \mathbf{D}_0$ are referred to as the \textit{conservative part} and \textit{dissipative part} of $\mathbf{F}_r$, respectively. By abuse of language, we also refer to $h$ and $\mu$ as the conservative and dissipative parts, respectively, and we will denote:
\[
\text{Cons}(\mathbf{F}_r) := h \quad \text{and} \quad \text{Diss}(\mathbf{F}_r) := \mu.
\]

\subsection{Monodromy and center problem.}

As it has been noted in the introduction, a singular point is called \textit{monodromic} when the orbits rotate around it. That is, in the case of analytical systems, a monodromic singular point is either a center or a focus. As we will see later, the factors of $h$, the function of the conservative part of $\mathbf{F}_r$, characterize the monodromy of the singular point. For this reason, we describe some aspects of the vector space of quasi-homogeneous polynomials and their factorizations in $\mathbb{C}[x, y]$.

Given a fixed $\mathbf{t} = (t_1, t_2)$, any quasi-homogeneous polynomial of type $\mathbf{t}$ and degree $k \in \mathbb{N}$, $p_k \in P^t_k$, can be expressed as:
\[
p_k(x, y) = x^{k_1} y^{k_2} \sum_{j=0}^{k_3} \alpha_j x^{t_2(k_3 - j)} y^{j t_1},
\]
where $k_1, k_2, k_3 \in \mathbb{N} \setminus \{0\}$, $k_1 < t_2$, $k_2 < t_1$, and $k = k_1 t_1 + k_2 t_2 + k_3 t_1 t_2$.

On the other hand, $p_k(x, y)$ is also associated with another homogeneous polynomial $p^{\text{hom}}_k$ of degree $k_3$ in the variables $X = x^{t_2}$ and $Y = y^{t_1}$, via the relation:
\[
p_k(x, y) := x^{k_1} y^{k_2} p^{\text{hom}}_k(X, Y),
\]
that is,
\[
p^{\text{hom}}_k(X, Y) := \frac{p_k(X^{1/t_1}, Y^{1/t_2})}{X^{k_1/t_2} Y^{k_2/t_1}}.
\]
We can get an expression of $h$ as a product of irreducible factors in $\mathbb{C}[x, y]$. Also, by scaling, we can assume that the leading coefficient of $h$ in the variable $y$ is one, that is:
\[
h(x, y) = x^{m_x} y^{m_y} \prod_{j=1}^{m} (y^{t_1} - \lambda_j x^{t_2})^{m_j}, \tag{3}
\]
where $m, m_x, m_y \in \mathbb{N} \cup \{0\}$; $m_j \in \mathbb{N}$ for $j = 1, \dots, m$; $\lambda_j \in \mathbb{C} \setminus \{0\}$ for $1 \leq j \leq m$, and $\lambda_i \neq \lambda_j$ for $i \neq j$ (for more details, see [26]).

\medskip

\textbf{Definition 1.} The polynomial $h \in P^t_k$ is \textit{monodromic} if the decomposition (3) of $h$ in $\mathbb{C}[x, y]$ verifies that $m_x$ and $m_y$ are null, $\lambda_j \in \mathbb{C} \setminus \mathbb{R}$, and $h$ is non-constant. That is, $h$ is non-constant and only has complex (no real) factors in its decomposition in $\mathbb{C}[x, y]$.

\medskip

The next result characterizes the monodromy of a system by means of the monodromy of the conservative part of its first quasi-homogeneous component.

\medskip

\textbf{Theorem 2.} Let $\mathbf{F}_r \in Q^t_r$. If $\text{Cons}(\mathbf{F}_r)$ is monodromic, then the origin of the system
\[
\dot{x} = \mathbf{F}(x) = \mathbf{F}_r(x) + \mathbf{F}_{r+1}(x) + \cdots, \quad \mathbf{F}_{r+j} \in Q^t_{r+j}, \, j \geq 1 \tag{4}
\]
is monodromic.

\begin{proof}
Using the conservative–dissipative decomposition of each of the quasi-homogeneous components of $\mathbf{F}$, this system can be written as:
\[
\dot{x} = \mathbf{F}(x) = \sum_{j=0}^{\infty} \left[ \mathbf{X}_{h_{r+j+|t|}} + \mu_{r+j} \mathbf{D}_0 \right], \tag{5}
\]
with $h_{r+j+|t|} \in P^t_{r+j+|t|}$ and $\mu_{r+j} \in P_{r+j}$.

We consider the change of variables:
\[
x = u^{t_1} \text{Cs}(\theta), \quad y = u^{t_2} \text{Sn}(\theta), \tag{6}
\]
where $(\text{Cs}(\theta), \text{Sn}(\theta))$ are the solutions of the initial value problem:
\[
\begin{cases}
\dot{x} = \mathbf{X}_H(x), \\
x(0) = (1, 0)^T, \quad \text{with } x = (x, y) \text{ and } H(x, y) = y^{2t_1} + x^{2t_2}.
\end{cases}
\]
We note that the point $(1, 0)^T$ belongs to the periodic ring, as $\mathbf{X}_H$ is a quasi-homogeneous global center. Therefore, $\text{Cs}(\theta)$ and $\text{Sn}(\theta)$ are periodic functions of period $T_t$.

Differentiating with respect to time, we obtain:
\[
\dot{x} = \frac{1}{u} \mathbf{D}_0 \dot{u} + \frac{1}{u^{2t_1 t_2 - |t|}} \mathbf{X}_H \dot{\theta}.
\]
Thus, 
\[
\dot{x} \wedge \mathbf{X}_H = \frac{1}{u} \mathbf{D}_0 \wedge \mathbf{X}_H \dot{u}, \quad \text{and} \quad \mathbf{D}_0 \wedge \dot{x} = \frac{1}{u^{2t_1 t_2 - |t|}} \mathbf{D}_0 \wedge \mathbf{X}_H \dot{\theta},
\]
where $\mathbf{D}_0 \wedge \mathbf{X}_H = 2 t_1 t_2 u^{2 t_1 t_2} \neq 0$.

In addition, using the conservative–dissipative decomposition of each quasi-homogeneous term, we get:
\[
\dot{x} \wedge \mathbf{X}_H = \sum_{j \geq 0} \mathbf{F}_{r+j} \wedge \mathbf{X}_H=
\] 
\[= \sum_{j \geq 0} \left[ \mathbf{X}_{h_{r+|t|+j}} + \mu_{r+j} \mathbf{D}_0 \right] \wedge \mathbf{X}_H
= \sum_{j \geq 0} \mathbf{X}_{h_{r+j+|t|}} \wedge \mathbf{X}_H + 2 t_1 t_2 H(x, y) \sum_{j \geq 0} \mu_{r+j}(x, y).
\]
and, 
\[
\mathbf{D}_0 \wedge \dot{x} =\sum_{j \geq 0} \mathbf{D}_0 \wedge \mathbf{X}_{h_{r+j+|t|}}(x, y) =
\] 
\[= \sum_{j \geq 0} (r + j + |t|) h_{r+j+|t|}(x, y).
\]

Moreover, for each $j \geq 0$, we have:
\[
h_{r+j+|t|}(x, y) = u^{r+j+|t|} h_{r+j+|t|}(\text{Cs}(\theta), \text{Sn}(\theta)) := u^{r+j+|t|} h_{r+j+|t|}(\theta),
\]
\[
\mu_{r+j}(x, y) = u^{r+j} \mu_{r+j}(\text{Cs}(\theta), \text{Sn}(\theta)) := u^{r+j} \mu_{r+j}(\theta),
\]
\[
\mathbf{X}_{h_{r+j+|t|}} \wedge \mathbf{X}_H = u^{r+j+2t_1 t_2 - |t|} \left[
-\frac{\partial h_{r+j+|t|}(\text{Sn}(\theta), \text{Cs}(\theta))}{\partial \text{Sn}(\theta)} \frac{d \text{Sn}(\theta)}{d \theta}
- \frac{\partial h_{r+j+|t|}(\text{Cs}(\theta), \text{Sn}(\theta))}{\partial \text{Cs}(\theta)} \frac{d \text{Cs}(\theta)}{d \theta}
\right] =\] 
\[=-\ u^{r+j+2t_1 t_2 - |t|} h'_{r+j+|t|}(\theta).
\]

Denoting $\text{Cons}(\mathbf{F}_r) = h = h_{r+|t|}$ and $\text{Diss}(\mathbf{F}_r) = \mu = \mu_r$, and applying the reparametrization of time $dt = \frac{2t_1 t_2}{u^r} d\tau$, the system (5) transforms into:
\[
u' = [2t_1 t_2 \mu(\theta) - h'(\theta)] u + O(u^2), \quad \theta' = (r + |t|) h(\theta) + O(u), \tag{7}
\]
where $u > 0$ and $' = \frac{d}{d\tau}$.
\end{proof}

Since \( h \) is monodromic, then \( h(\theta) \neq 0 \) for all \( \theta \in [0, T_t] \). Therefore, \( \theta' \neq 0 \) for all \( |u| \ll 1 \); consequently, the origin is a monodromic singular point of system (4).

\medskip

From Theorem 2, it is easy to prove the following:

\medskip

\textbf{Corollary 1.} \( h \) is monodromic if and only if \( \mathbf{X}_h \) has a center at the origin.

\medskip

The next statement establishes a sufficient and necessary condition for the monodromy of a quasi-homogeneous vector field.

\medskip

\textbf{Theorem 3.} Let \( \mathbf{F}_r \in Q^t_r \). The origin of \( \mathbf{F}_r \) is monodromic if and only if \( \mathbf{X}_h \) has a center at the origin, where \( h = \text{Cons}(\mathbf{F}_r) \).

\begin{proof}
From Lemma 1, we have \( \mathbf{F}_r = \mathbf{X}_h + \mu \mathbf{D}_0 \) with \( h \in P^t_{r+|t|} \) and \( \mu \in P^t_r \). The change of variables (6) and the reparametrization in time \( dt = \frac{2t_1 t_2}{u^r} d\tau \) transform the system \( \dot{x} = \mathbf{F}_r(x) \) into:
\[
u' = [2t_1 t_2 \mu(\theta) - h'(\theta)] u,
\]
\[
\theta' = (r + |t|) h(\theta), \tag{8}
\]
with \( u > 0 \) and \( ' = \frac{d}{d\tau} \).

To prove the necessary condition, we assume that \( \mathbf{X}_h \) does not have a center at the origin. From Corollary 1, \( h \) is not monodromic, which implies that \( h \) has a real factor in its decomposition in \( \mathbb{C}[x, y] \). Therefore, there exists a \( \theta_0 \in [0, T) \) such that \( h(\theta_0) = 0 \), meaning \( \theta = \theta_0 \) is invariant for (8). Consequently, the real factor is a solution of \( \dot{x} = \mathbf{F}_r(x) \), which leads to the conclusion that the origin of \( \mathbf{F}_r \) is non-monodromic, resulting in a contradiction.

The sufficient condition follows from Theorem 2 and Corollary 1.

\end{proof}

\medskip

From Theorem 3.7 of [27] and Theorem 3.3 of [26], we can obtain the following result that we will use in the next section.

\medskip

\textbf{Theorem 4.} Assume that \( \mathbf{F}_r = \mathbf{X}_h + \mu \mathbf{D}_0 \), with \( h \in P^t_r \) being monodromic. Then the origin of (1) satisfies:
\begin{itemize}
    \item[(a)] a global center if and only if \( f_0 = 0 \).
    \item[(b)] an unstable focus if and only if \( hf_0 > 0 \).
    \item[(c)] a stable focus if and only if \( hf_0 < 0 \).
\end{itemize}
where
\[
f_0 = -2\pi \sum_{{\text{Im}(\lambda_j) > 0}, \ h_{\text{hom}}(1, \lambda_j)=0} \text{Im} \left( \text{Res} \left[ \frac{\mu_{\text{hom}}(1, y)}{h_{\text{hom}}(1, y)}, \lambda_j \right] \right),
\]
and \( h_{\text{hom}}(x, y) \) and \( \mu_{\text{hom}}(x, y) \) are defined in (2).

The following theorem provides a necessary condition for a perturbed vector field of \( \mathbf{F}_r \) to have a center at the origin.

\textbf{Theorem 5.} If the origin of \( \dot{x} = \mathbf{F}_r(x) \) is a focus, then the origin of system (4) is also a focus with the same stability.

\begin{proof}
The orbits of system (7) in generalized polar coordinates are defined by the generalized Abel equation:
\[
\frac{du}{d\theta} =
\begin{bmatrix}
\frac{2t_1 t_2 \mu(\theta)}{(r + |t|) h(\theta)} - \frac{h'(\theta)}{(r + |t|) h(\theta)} \\
\end{bmatrix} u + O(u^2), \tag{9}
\]
with \( h(\theta) \neq 0 \) for all \( \theta \in [0, T_t] \) because it is monodromic.

If we denote \( u(\theta, u_0) = \sum_{n \geq 1} a_n(\theta) u_0^n \), the solution of Eq. (9) satisfying \( u(0, u_0) = u_0 \) can be found by replacing \( u(\theta, u_0) \) in Eq. (9), leading to:
\[
a_1(\theta) = \exp \left( \int_0^{\theta} \left[ \frac{2t_1 t_2 \mu(\alpha)}{(r + |t|) h(\alpha)} - \frac{h'(\alpha)}{(r + |t|) h(\alpha)} \right] d\alpha \right).
\]

The Poincaré map of system (4) is given by:
\[
P(u_0) = u(T_t, u_0) = \prod_{n \geq 1} a_n(T_t) u_0^n,
\]
defined for \( u_0 > 0 \). For system (1), the Poincaré map is given by \( P(u_0) = a_1(T_t) u_0 \), where
\[
a_1(T_t) = e^{-\frac{1}{r + |t|} \int_0^{T_t} \frac{h'(\theta)}{h(\theta)} d\theta + \frac{2t_1 t_2}{r + |t|} \int_0^{T_t} \frac{\mu(\theta)}{h(\theta)} d\theta} = e^{\frac{2t_1 t_2}{r + |t|} f_0}, \tag{10}
\]
and \( f_0 = \int_0^{T_t} \frac{\mu(\theta)}{h(\theta)} d\theta \).

Therefore, the first Lyapunov constant is the same for both systems (1) and (4). In the case that the origin of \( \dot{x} = \mathbf{F}_r(x) \) is a focus, we have \( a_1(T_t) \neq 1 \); otherwise, the origin is a center. Thus, the origin of the perturbed system is also a focus with the same stability as the focus of the non-perturbed system.
\end{proof}

\medskip

\textbf{Remark:} This result shows that if the origin of system (4) is a center, then the origin of system (1) must also be a center. This has been a motivation for studying the centers of vector fields of the form (1).\\

At first we give some necessary definitions.

\begin{itemize}
    \item An \textbf{involution} is a function \( \sigma \in C^\omega(U_0 \subset \mathbb{R}^2, \mathbb{R}^2) \), such that \( \sigma \circ \sigma = \text{Id} \), where \( U_0 \) is a neighborhood of the origin.
    \item A system \( \dot{x} = \mathbf{F}(x) \) (or a vector field \( \mathbf{F} \)) is \textbf{reversible} if there exists an involution \( \sigma \) (with \( \sigma \neq \pm \text{Id} \), \( \sigma(0) = 0 \)) such that \( \sigma^*\mathbf{F} = -\mathbf{F} \), where \( \sigma^*\mathbf{F} \) denotes the pull-back of \( \mathbf{F} \) by the transformation \( \sigma \).
    \item A system \( \dot{x} = \mathbf{F}(x) \) (or a vector field \( \mathbf{F} \)) is \textbf{orbitally reversible} if there exists an analytical scalar function \( f \), \( f(0) = 1 \), such that \( f \mathbf{F} \) is a reversible vector field.
    \item For quasi-homogeneous vector fields, the concepts of reversibility and orbital reversibility agree.
    \item A system \( \dot{x} = \mathbf{F}(x) \) (or a vector field \( \mathbf{F} \)) is \textbf{axis-reversible} if it is reversible to the involution \( \sigma(x, y) = (-x, y) \) (Rx-reversible) or to the involution \( \sigma(x, y) = (x, -y) \) (Ry-reversible).
\end{itemize}
The following proposition and corollary can be deduced from the works developed by Montgomery and Bochner (see [28]). For the sake of completeness, we give the proofs.\\

\textbf{Proposition 6.} Let \( \sigma \) be an involution, \( \sigma \in C^\omega \) and \( \sigma \neq \pm \text{Id} \). Then there exists \( \Psi(x, y) = (x + by + \cdots, ax + y + \cdots) \) such that \( \Psi \circ \sigma \circ \Psi^{-1}(x, y) \) is either \( (-x, y) \) or \( (x, -y) \).

\begin{proof}
We consider \( \sigma(x, y) = Ax + \cdots \). Since \( \sigma \) is an involution, \( A^2 = I_2 \). Using the linear transformation \( \phi(x, y) = (x + by, ax + y) \), this involution can be expressed in the new variables as:
\begin{enumerate}
    \item \( \sigma(x, y) = \pm (x + f(x, y), y + g(x, y)) \),
    \item \( \sigma(x, y) = (x + f(x, y), -y + g(x, y)) \),
    \item \( \sigma(x, y) = (-x + f(x, y), y + g(x, y)) \).
\end{enumerate}

The only involutions of type 1 are \( \pm \text{Id} \). Consider \( f(x, y) = A_{20}x^2 + A_{11}xy + A_{02}y^2 + \cdots \) and \( g(x, y) = B_{20}x^2 + B_{11}xy + B_{02}y^2 + \cdots \). Therefore,
\[
\sigma(x, y) = \left( x + A_{20}x^2 + A_{11}xy + A_{02}y^2 + \cdots, y + B_{20}x^2 + B_{11}xy + B_{02}y^2 + \cdots \right),
\]
imposing \( \sigma^2(x, y) = (x, y) \) it is easy to deduce that \( f \equiv 0 \equiv g \); consequently, the involutions that correspond to type 1 are \( \pm \text{Id} \).

For the involutions of type 2, we consider the change of variables \( \Psi(x, y) = \left( \frac{x + \sigma_1(x, y)}{2}, \frac{y - \sigma_2(x, y)}{2} \right) = (u, v) \) with \( \sigma = (\sigma_1, \sigma_2) \). Then,
\[
(\Psi \circ \sigma \circ \Psi^{-1})(u, v) = \Psi \circ \sigma(x, y) = \Psi(\sigma_1(x, y), \sigma_2(x, y)) = \left( \frac{x + \sigma_1(x, y)}{2}, \frac{-y + \sigma_2(x, y)}{2} \right) = (u, -v).
\]

The process is similar for the involutions of type 3.
\end{proof}

From Proposition 6, it is easy to prove the following:\\

\textbf{Corollary 2.} Let \( \mathbf{F} \) be a reversible vector field. Then there exists \( \Psi(x, y) = (x + by + \cdots, ax + y + \cdots) \in C^\omega \) such that \( \Psi^* \mathbf{F} \) is axis-reversible.\\

Before showing the following theorem, we need a technical lemma:\\

\textbf{Lemma 7.} Let \( t \neq (1, 1) \). The change of variables \( \Psi(x, y) = (x + by + \cdots, ax + y + \cdots) \in C^\omega \) can be decomposed in the form
\[
\Psi = (id + \tilde{\Psi}_{\leq 0}) \circ (id + \tilde{\Psi}_{> 0}),
\]
where \( \tilde{\Psi}_{\leq 0} = \sum_{i=0}^{-M} \tilde{\Psi}_i \), \( \tilde{\Psi}_{> 0} = \sum_{i=1}^{\infty} \tilde{\Psi}_i \), with \( \tilde{\Psi}_i \in \mathcal{Q}_t^i \), and \( \text{diag}(D \tilde{\Psi}_0(0)) = 0 \).

\begin{proof}
It is clear that \( \Psi \) can be expressed as
\[
\Psi = \tilde{\Psi}_{\leq 0} + id + \tilde{\Psi}_{> 0},
\]
with \( \text{diag}(D\Psi_0(0)) = 0 \). To prove the decomposition of \( \Psi \), it is sufficient to consider the following change of variables recursively defined:
\[
\Psi(M) := \Psi, \quad \Psi(i) := (id + \Psi(i+1)^{-i-1})^{-1} \circ \Psi(i+1) \quad \text{for } i = M - 1, \ldots, 0,
\]
i.e., we remove degree by degree the terms with negative or zero degree from \( \Psi \). It has been verified that \( -i \leq \min(\deg(\Psi(i))) \) for \( i = M, \ldots, 0 \).

Therefore, if we define
\[
\tilde{\Psi}_{> 0} := (id + \Psi(0)_{0})^{-1} \circ \Psi(0) - id
\]
and
\[
\tilde{\Psi}_{\leq 0} := (id + \Psi(0)_{0}) \circ \cdots \circ (id + \Psi(M)_{M}) - id,
\]
we obtain
\[
\Psi = (id + \tilde{\Psi}_{\leq 0}) \circ (id + \tilde{\Psi}_{> 0}),
\]
with \( 0 < \min(\deg(\tilde{\Psi}_{> 0})) \), \( \max(\deg(\tilde{\Psi}_{\leq 0})) \leq 0 \), and \( \text{diag}(D \tilde{\Psi}_{\leq 0}(0)) = 0 \).

Note that the change of variables with negative degree and their inverse, for the case \( t = (1, n) \), \( n > 1 \), are
\[
id + \tilde{\Psi}_{\leq 0} = \begin{cases}
u = x, \\
v = \alpha_1 x + \alpha_2 x^2 + \cdots + \alpha_n x^n + y,
\end{cases}
\]
and
\[
(id + \tilde{\Psi}_{\leq 0})^{-1} = id - \tilde{\Psi}_{\leq 0} = \begin{cases}
u = x, \\
v = -\alpha_1 x - \alpha_2 x^2 - \cdots - \alpha_n x^n + y.
\end{cases}
\]

For the change of variables with negative degree, for the case \( t = (m, n) \), \( 1 < m < n \), they are
\[
id + \tilde{\Psi}_{\leq 0} = \begin{cases}
u = x, \\
v = \alpha_1 x + \alpha_2 x^2 + \cdots + \alpha_M x^M + y,
\end{cases}
\]
with \( M = \binom{n}{m} \) and \( \lfloor \cdot \rfloor \) denoting the integer part.
\end{proof}

The next statement provides a sufficient and necessary condition about the reversibility of a quasi-homogeneous vector field.\\

\textbf{Theorem 8.} Let \( F_r \in \mathcal{Q}^t_r \). Then \( F_r \) is reversible if and only if there exists \( \Psi = Id + \Psi_0 \), where \( \Psi_0 \in \mathcal{Q}^t_0 \) and \( \text{diag}(D\Psi_0(0)) = 0 \), such that \( \Psi^* F_r \) is axis-reversible.

\begin{proof}
The sufficient condition is clear. We now prove the necessary condition depending on the types \( t \):

1. If \( t = (1, 1) \), we are in the case described in Corollary 2.

2. If \( t \neq (1, 1) \), from Corollary 2 and Lemma 7, there exists a change of variables \( (id + \Psi_{\leq 0}) \) and \( (id + \Psi_{> 0}) \) such that
\[
(id + \Psi_{\leq 0}) \circ (id + \Psi_{> 0})^* F_r = G,
\]
where \( G \) is axis-reversible. Considering the degrees of both vector fields, we have
\[
G = G_s + \cdots + G_r + \cdots,
\]
with \( s < r \) and \( G_j \in \mathcal{Q}^t_j \). Therefore,
\[
(id + \Psi_{> 0})^* F_r = (id - \Psi_{\leq 0})^* G.
\]
By decomposing it into quasi-homogeneous components, we find
\[
F_r + \tilde{F}_{r+1} + \cdots = \tilde{G}_l + \cdots + \tilde{G}_s + \cdots + \tilde{G}_r + \cdots,
\]
which leads to the conclusion that
\[
\tilde{G}_l = \cdots = \tilde{G}_s = \cdots = \tilde{G}_{r-1} \equiv 0.
\]
Moreover, since \( (id - \Psi_{\leq 0})^* G = \tilde{G}_r + \cdots \), we have
\[
\Psi^{-1} = \cdots = \Psi_{1-n} \equiv 0.
\]
Consequently, we find that \( id + \Psi_{\leq 0} = id + \Psi_0 \) and
\[
(id + \Psi_0) \circ (id + \Psi_{> 0})^* F_r = G,
\]
which implies that \( G \) is axis-reversible. Thus, in particular, \( (id + \Psi_0)^* F_r \) is axis-reversible.
\end{proof}

The following proposition establishes the axis-reversibility of a quasi-homogeneous vector field in function of its
conservative and dissipative part.\\

\textbf{Proposition 9.}
Assume that \( F_r = X_h + \mu D_0 \in \mathcal{Q}^t_r \).

\begin{itemize}
    \item[(a)] \( F_r \) is \( R_x \)-reversible if and only if \( h(-x, y) = h(x, y) \) and \( \mu(-x, y) = -\mu(x, y) \).
    \item[(b)] \( F_r \) is \( R_y \)-reversible if and only if \( h(x, -y) = h(x, y) \) and \( \mu(x, -y) = -\mu(x, y) \).
\end{itemize}

\begin{proof}
We will prove item (a), as the proof for the other case is similar.

Given \( F_r = (P, Q)^T \), \( F_r \) is \( R_x \)-reversible if and only if 
\[
P(-x, y) = P(x, y) \quad \text{and} \quad Q(-x, y) = -Q(x, y).
\]
To derive this, we compute:
\[
h(-x, y) = \frac{1}{r + |t|} D_0(-x, y) \wedge F_r(-x, y) = \frac{1}{r + |t|} \left( -t_1 x Q(-x, y) - t_2 y P(-x, y) \right).
\]
Thus, we can express:
\[
h(-x, y) = \frac{1}{r + |t|} \left( -t_1 x Q(-x, y) - t_2 y P(-x, y) \right) = h(x, y).
\]
Next, we analyze \( \mu \):
\[
\mu(-x, y) = \frac{1}{r + |t|} \left( \frac{\partial P(-x, y)}{\partial (-x)} + \frac{\partial Q(-x, y)}{\partial y} \right).
\]
This simplifies to:
\[
\mu(-x, y) = \frac{1}{r + |t|} \left( -\frac{\partial P(x, y)}{\partial x} - \frac{\partial Q(x, y)}{\partial y} \right) = -\mu(x, y).
\]

This shows the necessary condition for \( F_r \) being \( R_x \)-reversible.

Conversely, if \( h(-x, y) = h(x, y) \) and \( \mu(-x, y) = -\mu(x, y) \), we have:
\[
P(-x, y) = -\frac{\partial h(-x, y)}{\partial y} + t_1(-x) \mu(-x, y) = \frac{\partial h(x, y)}{\partial y} + t_1 x \mu(x, y) = P(x, y).
\]
Similarly,
\[
Q(-x, y) = \frac{\partial h(-x, y)}{\partial (-x)} + t_2 y \mu(-x, y) = -\frac{\partial h(x, y)}{\partial x} - t_2 y \mu(x, y) = -Q(x, y).
\]

This completes the proof.
\end{proof}
\section{Applications to the Degenerate Center Problem.}

In this section, we study the center set (center problem) of the vector field \( F_r \in \mathcal{Q}^t_r \) for \( 0 \leq r \leq 4 \), which involves determining when the origin of \( F_r \) is a center. It is known that if \( F_r \) is monodromic (i.e., the origin of \( F_r \) is monodromic) and either it is reversible or analytically integrable, then \( F_r \) is a center (i.e., the origin of \( F_r \) is a center).

In this section, we also study the subsets of centers which are either reversible or analytically integrable. The analysis of these systems will be conducted by classifying them by degrees, and for each degree, we will consider the expression of \( F_r \) depending on the type  $t$.\\

({\bf Case r = 0}) The system $\dot{\x}=\F_0(\x)$,
depending on the type $\t$, can be expressed as\\

1. For \( t = (1, 1) \):\[
\begin{cases}
\dot{x} = a_{10}x + a_{01}y, \\
\dot{y} = b_{10}x + b_{01}y,
\end{cases}
\tag{11}
\]

2. For \( t = (1, t_2) \):
\[
\begin{cases}
\dot{x} = a_{10}x, \\
\dot{y} = b_{t_2 0} x^{t_2} + b_{01}y, \quad t_2 > 1, \; t_2 \in \mathbb{N}.
\end{cases}
\tag{12}
\]

The following theorem solves the center, reversibility, and analytical integrability problems for systems (11) and (12).\\

\textbf{Theorem 10.}

1. The origin of system (12) is non-monodromic.

2. The origin of system (11) is monodromic if and only if 
\[
(b_{01} - a_{10})^2 + 4b_{10}a_{01} < 0,
\]
and in this case, it is a center if and only if 
\[
a_{10} = -b_{01}.
\]

3. If the origin of system (11) is a center, then system (11) is reversible and analytically integrable.\\

\begin{proof}
The Hamiltonian part of the conservative–dissipative decompositions of systems (11) and (12) are:

For system (11):
\[
h(x, y) = \frac{1}{2} \left[ b_{10}x^2 + (b_{01} - a_{10})xy - a_{01}y^2 \right].
\]

For system (12):
\[
h(x, y) = \frac{1}{t_2 + 1} \left[ b_{t_2 0}x^{t_2} + (b_{01} - a_{10}t_2)y \right].
\]

1. From Theorem 2, the origin of system (12) is non-monodromic because \( h \) has the real factor \( x \).

2. From Theorem 2, imposing the monodromy condition on system (11), we obtain the relation 

\[
(b_{01} - a_{10})^2 + 4a_{01}b_{10} < 0.
\]

On the other hand, if \( h(x, y) \) is monodromic, then it has the form 

\[
h(x, y) = [(y - ax)^2 + b^2x^2], \quad b \neq 0.
\]

From Theorem 4, we obtain 

\[
f_0 = \frac{\pi(a_{10} + b_{01})}{b}.
\]

Therefore, the origin of system (11) is a center if and only if 

\[
a_{10} = -b_{01}.
\]

3. In this case, system (11) is a Hamiltonian vector field; consequently, it is analytically integrable.

Finally, by Theorem 8, the change of variables 

\[
u = x, \quad v = -ax + y 
\]

transforms the system (11) into an axis-reversible system.
\end{proof}

\textbf{(Case \( r = 1 \))} The system $\bf{\dot{x}} = F_1(x)$, depending on the type $\t$, can be expressed as follows:\\

1. For \( t = (1, 1) \):
\[
\begin{cases}
\dot{x} = a_{20}x^2 + a_{11}xy + a_{02}y^2, \\
\dot{y} = b_{20}x^2 + b_{11}xy + b_{02}y^2,
\end{cases}
\tag{13}
\]

2. For \( t = (2, 3) \):
\[
\begin{cases}
\dot{x} = a_{01}y, \\
\dot{y} = b_{20}x^2,
\end{cases}
\tag{14}
\]

3. For \( t = (1, 2) \):
\[
\begin{cases}
\dot{x} = a_{20}x^2 + a_{01}y, \\
\dot{y} = b_{30}x^3 + b_{11}xy.
\end{cases}
\tag{15}
\]

The following result deals with the center problem for systems (13)–(15).\\

\textbf{Theorem 11}

1. The origin of systems (13) and (14) is non-monodromic. The origin of system (15) is monodromic if and only if 

\[
(b_{11} - 2a_{20})^2 + 8b_{30}a_{01} < 0.
\]

2. If the origin of system (15) is monodromic, then it is a center if and only if 

\[
2a_{20} = -b_{11}.
\]

3. If the origin of system (15) is a center, then system (15) is reversible and analytically integrable.

\begin{proof}
Systems (13)–(15) can be expressed in the form \( Xh + \mu D_0 \) with:

For system (13):
\[
\begin{cases}
h(x, y) = \frac{1}{3} \left[ b_{20}x^3 + (b_{11} - a_{20})x^2y + (b_{02} - a_{11})xy^2 - a_{02}y^3 \right], \\
\mu(x, y) = \frac{1}{3} \left[ x(2a_{20} + b_{11}) + y(a_{11} + 2b_{02}) \right].
\end{cases}
\]

For system (14):
\[
\begin{cases}
h(x, y) = \frac{1}{6} \left( 2b_{20}x^3 - 3a_{01}xy \right), \\
\mu(x, y) = 0.
\end{cases}
\]

For system (15):
\[
\begin{cases}
h(x, y) = \frac{1}{4} \left[ b_{30}x^4 + (b_{11} - 2a_{20})x^2y - 2a_{01}y^2 \right], \\
\mu(x, y) = \frac{1}{4} \left[ (2a_{20} + b_{11})x \right].
\end{cases}
\]

1. The Hamilton function of systems (13) and (14) is a polynomial of degree three and it has a real factor in its decomposition (3). Therefore, from Theorem 2, the origin of system (13) is non-monodromic. On the other hand, the Hamilton function of the decomposition of system (15) has only complex factors if and only if 

\[
(b_{11} - 2a_{20})^2 + 8b_{30}a_{01} < 0.
\]

2. Assume that system (15) is monodromic. Then, the Hamiltonian part of the decomposition of this system can be expressed in the form 

\[
h(x, y) = [(y - ax^2)^2 + b^2x^4], \quad b \neq 0,
\]

and using Theorem 4, we obtain 

\[
f_0 = -\frac{(2a_{20} + b_{11})}{2b}.
\]

Therefore, the origin of system (15) is a center if and only if 

\[
2a_{20} + b_{11} = 0.
\]

3. In the case of system (15), if it has a center at the origin, then it is Hamiltonian and consequently analytically integrable. On the other hand, the change of variables 

\[
u = x, \quad v = -ax^2 + y 
\]

transforms the system (15) into an axis-reversible system. From Theorem 8, the system is reversible.
\end{proof}

\textbf{(Case \( r = 2 \))} The system \(\dot{x} = F_2(x)\), depending on the type \( t \), can be expressed as follows:\\

1. For \( t = (1, 1) \):
\[
\begin{cases}
\dot{x} = a_{30}x^3 + a_{21}x^2y + a_{12}xy^2 + a_{03}y^3, \\
\dot{y} = b_{30}x^3 + b_{21}x^2y + b_{12}xy^2 + b_{03}y^3,
\end{cases}
\tag{16}
\]

2. For \( t = (1, 2) \):
\[
\begin{cases}
\dot{x} = a_{30}x^3 + a_{11}xy, \\
\dot{y} = b_{40}x^4 + b_{21}x^2y + b_{02}y^2,
\end{cases}
\tag{17}
\]

3. For \( t = (1, 3) \):
\[
\begin{cases}
\dot{x} = a_{30}x^3 + a_{01}y, \\
\dot{y} = b_{50}x^5 + b_{21}x^2y.
\end{cases}
\tag{18}
\]

To study the family (16), we reduce the number of parameters for simplicity.\\

\textbf{Proposition 12} The origin of system (16) is monodromic if and only if there exists a degree zero change of variables and a time reparametrization that transforms it into

\[
\begin{cases}
\dot{x} = \mu_0x^3 + (\mu_1 - 2(B^2 + 1))x^2y + \mu_2xy^2 - 4y^3, \\
\dot{y} = 4B^2x^3 + \mu_0x^2y + (\mu_1 + 2(B^2 + 1))xy^2 + \mu_2y^3.
\end{cases}
\]

\begin{proof}
The functions \( h \) and \( \mu \) of the conservative and dissipative parts of the decomposition of system (16) are given by:

\[
\begin{cases}
h(x, y) = \frac{1}{4} \left[ b_{30}x^4 + (b_{21} - a_{30})x^3y + (b_{12} - a_{21})x^2y^2 + (b_{03} - a_{12})xy^3 - a_{03}y^4 \right], \\
\mu(x, y) = \frac{1}{4} \left[ (3a_{30} + b_{21})x^2 + (2a_{21} + 2b_{12})xy + (a_{12} + 3b_{03})y^2 \right].
\end{cases}
\]

From Theorem 2, the origin of (16) is monodromic if and only if \( h \) is monodromic. In this case, we can assume that \( a_{03} \neq 0 \); otherwise, \( h \) is non-monodromic. By scaling in time, we can assume the coefficient of \( y^4 \) in the polynomial \( h \) is 1. 

We will distinguish two different possibilities:

(a) If \( h \) has simple factors in \( \mathbb{C}[x, y] \), then \( h \) has the form 

\[
h(x, y) = [(y - a_1x)^2 + b_1^2x^2][(y - a_2x)^2 + b_2^2x^2], \quad b_i \neq 0, \quad i = 1, 2, \quad (a_1, b_1) \neq (a_2, b_2).
\]

The change of variables 

\[
u = (|b_1| - a_1A)x + Ay, \quad v = y - (a_1 + A|b_1|)x,
\]

where \( A = 0 \) if \( a_1 = a_2 \) or \( A \) satisfies the equation 

\[
\frac{a_2 - a_1}{b_1}A^2 + \left( 1 - \frac{a_2 - a_1}{b_1^2} + \frac{b_2^2}{b_1^2} \right)A - \frac{a_2 - a_1}{b_1} = 0
\]

in the other case, transforms system (16) into 

\[
\dot{x} = X\tilde{h} + \tilde{\mu} D_0,
\]

where 

\[
\tilde{h}(x, y) = (y^2 + x^2)(y^2 + B^2x^2),
\]

\[
\tilde{\mu}(x, y) = \mu_0x^2 + \mu_1xy + \mu_2y^2.
\]

Writing the system in the usual form yields (19).

(b) If \( h \) has multiple factors in \( \mathbb{C}[x, y] \), then \( h \) can be written as 

\[
h(x, y) = [(y - a_1x)^2 + b_1^2x^2]^2.
\]

Using the change of variables 

\[
x_1 = |b_1|x, \quad x_2 = y - a_1x,
\]

\( h \) and \( \mu \) can be transformed into:

\[
\tilde{h}(x, y) = (y^2 + x^2)^2,
\]

\[
\tilde{\mu}(x, y) = \mu_0x^2 + \mu_1xy + \mu_2y^2,
\]

and writing the system in the usual form gives (19) in the particular case \( B = 1 \).

\end{proof}

\textbf{Theorem 13}

1. The origin of system (17) is non-monodromic.

2. The origin of system (18) is monodromic if and only if 

\[
(b_{21} - 3a_{30})^2 + 12a_{01}b_{50} < 0.
\]

In this case, the origin of (18) is a center if and only if 

\[
3a_{30} + b_{21} = 0.
\]

These centers are reversible and analytically integrable.

3. The origin of system (19) is a center if and only if 

\[
\mu_0 = -B\mu_2.
\]

These centers are reversible if either \( B = 1 \) or \( \mu_0 = \mu_2 = 0 \), and are analytically integrable if and only if \( B \neq 1 \) and \( \mu_0 = \mu_1 = \mu_2 = 0 \).\\

\textbf{Remark 1}
Notice that there exist systems in the family (19) which are non-reversible and non-integrable centers. Moreover, there exist systems which are reversible and non-integrable and vice versa. This is a new situation compared to that of non-degenerate and nilpotent centers.\\

In the proof of this theorem, we will use the following results. Their demonstrations can be found in [26].\\

\textbf{Proposition 14. }
If \( F_r \) is irreducible and \( h \) has multiple factors in decomposition (3), \( h \neq \text{cte} \), and \( \mu \neq 0 \), then \( F_r \) is non-analytically integrable.\\

\textbf{Theorem 15.} Assume that \( F_r \) is irreducible, \( h \) is monodromic, and \( h \) has more than two factors, all of them simple, in its decomposition (3). The system \( \dot{x} = F_r(x) \) has a first integral if and only if either \( \mu \equiv 0 \) or there exist \( n_x, n_y, n_i \), \( i = 1, \ldots, m \), non-negative numbers, not all zero, such that:

\[
\begin{cases}
\text{Res}[\eta_{\text{hom}}(X, 1), 0] = \frac{(n_x + 1)(r + |t|) - M_0}{t_2M_0}, & \text{if } \delta_x = 1, \\
\text{Res}[\eta_{\text{hom}}(1, Y), 0] = -\frac{(n_y + 1)(r + |t|) - M_0}{t_1M_0}, & \text{if } \delta_y = 1, \\
\text{Res}[\eta_{\text{hom}}(1, Y), \lambda_i] = -\frac{(n_i + 1)(r + |t|) - M_0}{M_0}, & i = 1, \ldots, m,
\end{cases}
\]

with 

\[
\eta_{\text{hom}}(X, Y) = \frac{\mu_{\text{hom}}(X, Y)}{X^{\delta_x} Y^{\delta_y} h_{\text{hom}}(X, Y)}
\]

and 

\[
M_0 = t_1(n_x + 1)\delta_x + t_2(n_y + 1)\delta_y + t_1t_2 \sum_{j=1}^{m}(n_j + 1).
\]

Moreover, in this case, a first integral of degree \( M_0 \) is given by:

\[
U = x^{(n_x + 1)\delta_x}y^{(n_y + 1)\delta_y} \prod_{i=1}^{m}(y t_1 - \lambda_i x t_2)^{n_i + 1}.
\]

\begin{proof} [Theorem 13.]
The functions \( h \) and \( \mu \) of the conservative and dissipative part of the decomposition of systems (17)–(19) are:\\

For system (17):
\[
\begin{cases}
h(x, y) = \frac{1}{5}[b_{40}x^5 + (b_{21} - 2a_{30})x^3y + (b_{02} - 2a_{11})xy^2], \\
\mu(x, y) = \frac{1}{5}[(3a_{30} + b_{21})x^2 + (a_{11} + 2b_{02})y].
\end{cases}
\]

For system (18):
\[
\begin{cases}
h(x, y) = \frac{1}{6}[b_{50}x^6 + (b_{21} - 3a_{30})x^3y - 3a_{01}y^2], \\
\mu(x, y) = \frac{1}{6}(3a_{30} + b_{21})x^2.
\end{cases}
\]

For system (19):
\[
\begin{cases}
h(x, y) = (y^2 + x^2)(y^2 + B^2x^2), \\
\mu(x, y) = \mu_0x^2 + \mu_1xy + \mu_2y^2.
\end{cases}
\]

\begin{enumerate}
    \item The Hamiltonian function of the decomposition of system (17) has the real factor \( x \) in its factorization (3). Therefore, from Theorem 2, the origin of (17) is non-monodromic.

    \item From Theorem 2, imposing the monodromy condition on system (18), we obtain the relation:

    \[
    (b_{21} - 3a_{30})^2 + 12b_{50}a_{01} < 0.
    \]

    Assume that system (18) is monodromic. The Hamiltonian part of the decomposition of this system can be expressed in the form:

    \[
    h(x, y) = c[(y - ax^3)^2 + b^2x^6],
    \]

    where \( c = 3a_{01} \), \( a = \frac{b_{21} - 3a_{30}}{6a_{01}} \), and 

    \[
    b = \frac{\sqrt{|(b_{21} - 3a_{30})^2 + 12b_{50}a_{01}|}}{6a_{01}}.
    \]

    From Theorem 4, we obtain 

    \[
    f_0 = \frac{\pi b}{3a_{30} + b_{21}}.
    \]

    Thus, we can conclude that the origin of system (18) is a center if and only if 

    \[
    3a_{30} + b_{21} = 0.
    \]

    In the case where the origin of system (18) is a center, using Theorem 8, the change of variables 

    \[
    u = |b|x^3, \quad v = y - ax^3
    \]

transforms the system into an axis-reversible system. Moreover, these centers are integrable because the conditions for the center cancel out the dissipative part and the system is Hamiltonian.

    \item For studying system (19), we distinguish two subcases:

    \begin{enumerate}
        \item[(1)] \( B \neq 1 \). From Theorem 4, we obtain 

        \[
        f_0 = -\frac{\mu_0 + B\mu_2}{2B(B + 1)}.
        \]

        Therefore, the origin of system (19) is a center if and only if 

        \[
        \mu_0 + B\mu_2 = 0.
        \]
In the case that the origin of system (19) is a center, if \( \mu_0 = \mu_2 = 0 \) then the system (19) is axis-reversible; consequently, it is reversible. In the other case, i.e., \( \mu_0 + B\mu_2 = 0 \) and \( |\mu_0| + |\mu_2| \neq 0 \), a change of variables does not exist:

\[
\Psi(x, y) = (x + ay, bx + y)
\]

such that \( \Psi^*F_2 \) is axis-reversible; from Theorem 8, the system is non-reversible.

From Theorem 15, system (19) is integrable if the following conditions are verified:

\begin{enumerate}
    \item[(a)] 
    \[
    \text{Res}[\eta_{\text{hom}}(1, Y), i] = \frac{\mu_1 - (\mu_0 - \mu_2)i}{2(B^2 - 1)} = -\frac{4(n_1 + 1) - M_0}{M_0}.
    \]

    \item[(b)] 
    \[
    \text{Res}[\eta_{\text{hom}}(1, Y), -i] = \frac{\mu_1 + (\mu_0 - \mu_2)i}{2(B^2 - 1)} = -\frac{4(n_2 + 1) - M_0}{M_0}.
    \]

    \item[(c)] 
    \[
    \text{Res}[\eta_{\text{hom}}(1, Y), Bi] = \frac{\mu_1 B - (\mu_0 - B^2\mu_2)i}{2B(B^2 - 1)} = -\frac{4(n_3 + 1) - M_0}{M_0}.
    \]

    \item[(d)] 
    \[
    \text{Res}[\eta_{\text{hom}}(1, Y), -Bi] = \frac{\mu_1 B + (\mu_0 - B^2\mu_2)i}{2B(B^2 - 1)} = -\frac{4(n_4 + 1) - M_0}{M_0}.
    \]
\end{enumerate}

where 

\[
M_0 = (n_1 + 1) + (n_2 + 1) + (n_3 + 1) + (n_4 + 1).
\]

From here, we can conclude that the unique condition for integrability is \( \mu_0 = \mu_1 = \mu_2 = 0 \), i.e., when system (19) is a Hamiltonian system.

\item[(2)] \( B = 1 \)

From Theorem 4, we have:

\[
f_0 = -\frac{\pi}{2} (\mu_0 + \mu_2).
\]

Therefore, the origin of system (19) is a center if and only if 

\[
\mu_0 + \mu_2 = 0.
\]

In this case, to study the reversibility, we consider two subcases:

1. If \( \mu_0 = 0 \), then \( \mu_2 = 0 \) and, by applying Proposition 9, the system is axis-reversible.

2. If \( \mu_0 \neq 0 \), the change of variables 

\[
\text{id} + \Psi_0 = (x + \beta y, y - \beta x)
\]

with 

\[
\beta = \frac{\mu_1 \pm \sqrt{\mu_1^2 + 4\mu_0^2}}{2\mu_0}
\]

transforms the system into an axis-reversible system. From Theorem 8, it is reversible.

By applying Proposition 14, we can observe that \( h \) has double factors; consequently, the centers are non-integrable. 

\end{enumerate}
\end{enumerate}
\end{proof}
\textbf{Case (\( r = 3 \)).} The system \( \dot{x} = F_3(x) \), depending on the type \( t \), can be expressed as follows:

1. For \( t = (1, 1) \):
\[
\begin{cases}
\dot{x} = a_{40} x^4 + a_{31} x^3 y + a_{22} x^2 y^2 + a_{13} x y^3 + a_{04} y^4, \\
\dot{y} = b_{40} x^4 + b_{31} x^3 y + b_{22} x^2 y^2 + b_{13} x y^3 + b_{04} y^4.
\end{cases} 
\tag{21}
\]

2. For \( t = (1, 2) \):
\[
\begin{cases}
\dot{x} = a_{40} x^4 + a_{21} x^2 y + a_{02} y^2, \\
\dot{y} = b_{50} x^5 + b_{31} x^3 y + b_{12} x y^2.
\end{cases} 
\tag{22}
\]

3. For \( t = (1, 3) \):
\[
\begin{cases}
\dot{x} = a_{40} x^4 + a_{11} x y, \\
\dot{y} = b_{60} x^6 + b_{31} x^3 y + b_{02} y^2.
\end{cases} 
\tag{23}
\]

4. For \( t = (1, 4) \):
\[
\begin{cases}
\dot{x} = a_{40} x^4 + a_{01} y, \\
\dot{y} = b_{70} x^7 + b_{31} x^3 y.
\end{cases} 
\tag{24}
\]

5. For \( t = (2, 3) \):
\[
\begin{cases}
\dot{x} = a_{11} x y, \\
\dot{y} = b_{30} x^3 + b_{02} y^2.
\end{cases}
\tag{25}
\]

6. For \( t = (2, 5) \):
\[
\begin{cases}
\dot{x} = a_{01} y, \\
\dot{y} = b_{40} x^4.
\end{cases}
\tag{26}
\]

The following result analyzes the center problem for systems (22)–(26).\\

\textbf{Theorem 16}
\begin{enumerate}
    \item The origin of the families $(21)$–$(23)$, $(25)$ and $(26)$ is non-monodromic.
    \item The origin of system $(24)$ is monodromic if and only if 
    \[
    (b_{31} - 4a_{40})^2 + 16b_{70}a_{01} < 0.
    \]
    \item In the case that the origin of system $(24)$ is monodromic, it is a center if and only if 
    \[
    4a_{40} + b_{31} = 0.
    \]
    These centers are reversible and analytically integrable.
\end{enumerate}

\begin{proof}
\begin{enumerate}
    \item From Theorem 2, it is sufficient to verify that the Hamiltonian functions of the conservative–dissipative decomposition of these systems have a real factor.
    \item The functions $h$ and $\mu$ of the conservative–dissipative decomposition of system $(24)$ are:
    \[
    h(x, y) = \frac{1}{8} \left[ b_{70}x^8 + (b_{31} - 4a_{40})x^4y - 4a_{01}y^2 \right],
    \]
    \[
    \mu(x, y) = \frac{1}{8} (4a_{40} + b_{31})x^3.
    \]
    From Theorem 2, we conclude that system $(24)$ is monodromic if and only if
    \[
    (b_{31} - 4a_{40})^2 + 16b_{70}a_{01} < 0.
    \]
    \item Assume that the origin of system $(24)$ is monodromic. The Hamiltonian part of the decomposition of system $(24)$ can be expressed in the form
    \[
    h(x, y) = c \left[ (y - ax^4)^2 + b^2x^8 \right],
    \]
    where 
    \[
    c = 4a_{01}, \quad a = \frac{(b_{31}-4a_{40})^2}{8a_{01}}, \quad b = \frac{\sqrt{|(b_{31}-4a_{40})^2 - 16b_{70}a_{01}|}}{4a_{01}}.
    \]
    Using Theorem 4, we obtain 
    \[
    f_0 = \frac{\pi}{8} \left( \frac{4a_{40}+b_{31}}{b} \right),
    \]
    thus the origin of system $(24)$ is a center if and only if 
    \[
    4a_{40} + b_{31} = 0.
    \]
    In this case, system $(24)$ is Hamiltonian and, consequently, analytically integrable. Furthermore, from Theorem 8, the change of variables 
    \[
    u = |b|x^4, \quad v = y - ax^4
    \]
    transforms system $(24)$ into an axis-reversible system.
\end{enumerate}
\end{proof}
\textbf{(Case $r = 4$).} The system $\dot{x} = F_4(x)$ can be expressed for various types $t$: \\

1. For \( t = (1, 1) \):
 \[
 \begin{cases}
   \dot{x} = a_{50}x^5 + a_{41}x^4y + a_{32}x^3y^2 + a_{23}x^2y^3 + a_{14}xy^4 + a_{05}y^5, \\
   \dot{y} = b_{50}x^5 + b_{41}x^4y + b_{32}x^3y^2 + b_{23}x^2y^3 + b_{14}xy^4 + b_{05}y^5.
    \end{cases}
    \tag{27}
 \]

2. For t =(1, $t_2$), $t_2 > 1$:
 \[
 \begin{cases}
 \dot{x} = a_{50}x^5 + a_{31}\chi_{t_2=2}x^3y + a_{12}\chi_{t_2=2}xy^2 + a_{11}\chi_{t_2=4}xy + a_{01}\chi_{t_2=5}y, \\
 \dot{y} = b_{t_2+4,0}x^{t_2+4} + b_{41}x^4y + b_{22}\chi_{t_2=2}x^2y^2 + b_{03}\chi_{t_2=2}y^3 + b_{02}\chi_{t_2=4}y^2.
 \end{cases}
  \tag{28}
 \]
where \( \chi_{p} \) is equal to 1 if the proposition \( p \) is true, and is equal to 0 otherwise.\\

\textbf{Proposition 17.}  The origin of (27) is monodromic if and only if there exists a change of variables of degree zero and a scaling in
the time which transforms the system $\dot{\x} = Xh + \mu \D_0$ into a system with the form $\dot{\x} = X_{\tilde{h}} + \tilde{\mu} \D_0$, being $\tilde{h}$ and $\tilde{\mu}$ one and only
one of the following:

\[
\begin{cases}
\tilde{h}= (x^2 + y2)^3 \\
\tilde{\mu}= \mu_0x^4 + \mu_1x^3y +  \mu_2x^2y2 + \mu_3xy^3 +  \mu_4y^4
\end{cases}
\tag{29}
\]

\[
\begin{cases}
\tilde{h}= (x^2 + y^2)^2(x^2 + B^2y^2) \\
\tilde{\mu}= \mu_0x^4 + \mu_1x^3y +  \mu_2x^2y2 + \mu_3xy^3 +  \mu_4y^4
\end{cases}
\tag{30}
\]
 
\[
\begin{cases}
\tilde{h}= (x^2 + y^2)(x^2 + B^2y^2))[(y - Ax)^2 + C^2x^2], \\
\tilde{\mu}= \mu_0x^4 + \mu_1x^3y +  \mu_2x^2y2 + \mu_3xy^3 +  \mu_4y^4
\end{cases}
\tag{31}
\]
\begin{proof}
The function $h$ for system (27) is a homogeneous
polynomial of degree $6$, without loss of generality we can assume
that the coefficient of $y^6$ in $\hx$ is $1$, otherwise $\hx$ is non-monodromic. Possible options
are:
\begin{itemize}

\item $h$ has a pair of conjugate complex factors of multiplicity
3. In this case $h$ can be expressed in the form $h=[(y-a_1x)^2 +
b_1^2x^2]^3$. The change of variables $x=b_1x$, $v=y-a_1x$
transform system (27) into (29).

\item $h$ has a pair of conjugate complex factors of multiplicity
2 and a simple pair of  conjugate complex factor or  $h$ has three
simple pair of complex conjugate factors. By applying the changes
of variables described in the proof of item a) of Proposition
12, system (27) is transformed into
(30) and (31) respectively.
\end{itemize} 
\end{proof}

\textbf{Theorem 18.} 
\begin{enumerate}
\item The origin of system (28) is monodromic if and only
if $t_2$=5 and $(b_{41}-5a_{50})^2 + 20b_{90}a_{01}<0$. In this
case, it is a center if and only if $b_{41}+5a_{50}=0$. These
centers are reversible and analytically integrable.

\item The origin of system (29) is a center if and
only if \:$\mu_2=-3(\mu_4+\mu_0)$. In this case the system is
non-integrable. It is reversible if and only if any of the
following situations occurs:
\begin{itemize}
\item $\mu_0=\mu_2=\mu_4=0$ \item There exists
$\beta\in\mathbb{R}$, $\beta\neq0$ such that

\hspace{2cm} $\mu_2=\frac{-3}{2}\frac{\beta^2-1}{\beta}[\mu_1 +
\frac{\beta^2-1}{\beta}\mu_0]$

\hspace{2cm} $\mu_3=\frac{1}{2\beta^2}[\frac{(\beta^2-1)(\beta^4 -
6\beta^2+1)}{\beta}\mu_0 + (\beta^4-4\beta^2 + 1)\mu_1]$

\hspace{2cm}$\mu_4=\frac{1}{2\beta}[\frac{\beta^4-4\beta^2+1}{\beta}\mu_0
+ (\beta^2 - 1)\mu_1]$
\end{itemize}

\item The origin of system (30) is a center if and
only if $2\mu_4B^2+(\mu_2+\mu_4+\mu_0)B+2\mu_0=0$ with $B\neq0$.
In this case the system is non-integrable. It is reversible if and
only if $\mu_0=\mu_2=\mu_4=0$.

\item The origin of system (31) is a center if and
only if
$(C+B+{C}^{3}+{B}^{2}+3BC+2B{C}^{2}+2{C}^{2}+{B}^{2}C+{A}^{2}C)
\mu_0+ ( 2ABC+{B}^{2}A+AB)\mu_1+(B{C}^{3}+{A}^{2}BC+{B}^{2}{A}^{2}
+B{A}^{2}+{B}^{2}C+{B}^{2}{C}^{2}+B{C}^{2})\mu_2+ (
A{C}^{2}{B}^{2}+{A}^{3}{B}^{2}+ A{C}^{2}B+{A}^{3}B+2\,A{B}^{2}C
)\mu_3+ ( {A}^{4}B+{B}^
{3}C+3{B}^{2}{A}^{2}C+{A}^{2}BC+{B}^{2}{C}^{4}+2{B}^{3}{C}^{2}+{C}
^{4}B+2{B}^{2}{C}^{2}+B{C}^{3}+{B}^{3}{C}^{3}+{B}^{3}C{A}^{2}+2{B}
^{2}{A}^{2}{C}^{2}+3{B}^{2}{C}^{3}+{B}^{2}{A}^{4}+2{A}^{2}{C}^{2}B
 ) \mu_4=0$.
In this case the system is integrable if and only if
$\mu_0=\mu_1=\mu_2=\mu_3=\mu_4=\mu_5=0$, and is reversible if and
only if $A=0$ y $\mu_0=\mu_2=\mu_4=0$.

\end{enumerate}

\begin{proof}\hspace*{1cm}
\begin{enumerate}
\item The Hamiltonian function of the conservative-dissipative
decomposition of system (28) is

$(5+t_2)h=b_{t_2+4,0}x^{t_2+5}+(b_{41}-t_2a_{50})x^5y+\chi_{\{t_2=2\}}(b_{22}-2a_{31})x^3y^2+
\chi_{\{t_2=4\}}(b_{02}-4a_{11})xy^2-5a_{01}\chi_{\{t_2=5\}}y^2+\chi_{\{t_2=2\}}(b_{02}-2a_{12})xy^3.$\\
We distinguish two cases:

\begin{itemize}

\item If  $t_2\neq 5$, $h$ has the real factor $x$, from Theorem
2, system (28) is non-monodromic.

\item If $t_2$=5, $10h(x,y)=
b_{90}x^{10}+(b_{41}-5a_{50})yx^5-5a_{01}y^2$. In this case, from
Theorem 2, it is monodromic if and only if
$(b_{41}-5a_{50})^2 + 20b_{90}a_{01}<0$. In this case the
Hamiltonian part of decomposition of system (28) can be
expressed in the form $h=[(y-ax^5)^2 + b^2x^{10}], \:\: b\neq0$
then, from Theorem 4, we obtain $f_0=
\frac{\pi(b_{41}+5a_{50})}{b}$ and we can conclude that the origin
of system (28) is a center if and only if
$b_{41}+5a_{50}=0$. From Theorem 8, system
(28) is reversible because the change of variables $u=x$,
$v=y-ax^5$ transforms the system into an axis-reversible system.
In the case that the origin of system (28) is a center,
the system is Hamiltonian and consequently analytically
integrable.
\end{itemize}
\item  By applying Theorem 4, $f_0$ =
$\pi\frac{\mu_2+3\mu_4+3\mu_0}{8}$. Therefore, the origin of
system (29) is a center if and only if
$\mu_2+3\mu_4+3\mu_0=0$.

From Proposition 14 this system is non-integrable.

In this case, by applying Theorem 8, the only change
of variables that transforms the Hamilton part $h$ into an other
Hamiltonian part $\tilde{h}$ which is axis-reversible has the form
$u=x+\beta y$, $v=-\beta x + y$. This change transforms $h$ and
$\mu$ into

\begin{itemize}
\item $\tilde{h}(x,y)= (v^2+u^2)^3/(1+\beta^2)^2$

\item $\tilde{\mu}(x,y)=\frac{\mu_0\beta^4+\beta^3\mu_1+\mu_2\beta^2+\beta\mu_3+\mu_4}{(1+\beta^2)^4}v^4+\frac{\mu_4\beta^4-\beta^3\mu_3+\mu_2\beta^2-\beta\mu_1+\mu_0}{(1+\beta^2)^4}u^4$

    \hspace{1.5cm}$+\frac{\mu_2\beta^4+3\beta^3\mu_3-3\beta^3\mu_1-4\mu_2\beta^2+6\mu_0\beta^2+6\mu_4\beta^2+3\beta\mu_1-3\beta\mu_3+\mu_2}{(1+\beta^2)^4}v^2u^2$

    \hspace{1.5cm}$+\frac{-\beta^4\mu_3-4\mu_4\beta^3+2\mu_2\beta^3-3\beta^2\mu_1+3\beta^2\mu_3+4\mu_0\beta-2\mu_2\beta+\mu_1}{(1+\beta^2)^4}u^3v$

    \hspace{1.5cm}$+\frac{-\beta^4\mu_1+4\mu_0\beta^3-2\mu_2\beta^3-3\beta^2\mu_3+3\beta^2\mu_1-4\mu_4\beta+2\mu_2\beta+\mu_3}{(1+\beta^2)^4}uv^3$
\end{itemize}

and we can observe that $\tilde{h}$ is always even in $u$ and $v$.
Therefore, from Proposition 9, system
(29) is axis-reversible if there exists $\beta$ such
that $\tilde{\mu}$ is odd in $u$ or $v$, i.e. if there exist $\beta$ such that 
annul coefficients $u^4$, $v^4$ and $u^2v^2$ en $\tilde{\mu}$, this is the condition of the theorem.
\item From Theorem 4, we obtain $f_0=
\frac{\pi}{2}\frac{2\mu_4B^2+\mu_2B+\mu_4B+\mu_0B+2\mu_0}{(B+1)^2B}$.
Therefore, the origin of system (30) is a center if
and only if $2\mu_4B^2+(\mu_2+\mu_4+\mu_0)B+2\mu_0=0$.

In this case, system (30) is axis-reversible if and
only if $\mu_2=\mu_4=\mu_0=0$. Because, from Proposition 9.
\begin{itemize}

     \item $\tilde{h}=(x^2 + y^2)^2(x^2 + B^2y^2)$ is even in $x$ and $y$.
     \item $\tilde{\mu}=\mu_0x^4+\mu_1x^3y+\mu_2x^2y^2+\mu_3xy^3+\mu_4y^4$  is odd in $x$ or $y$ if and only if $\mu_2=\mu_4=\mu_0=0$

\end{itemize}

In other case, from Theorem 8, system
(30) is non-reversible because there is no change of
variables which transforms this system into an axis-reversible
system.

From Proposition 14 system (30) is
non-integrable.

\item From Theorem 4 we calculate the expression of $f_0$.
Imposing $f_0=0$ we obtain the result.

To study the integrability of system (31) we use
Theorem 15. We obtain $n_1$, $n_2$ y $n_3$ $\in$
$\mathbb{N}$, such that

 \hspace{2cm}Res[$\eta^{hom}(1,Y)$, $\pm i$] = $-\frac{6(n_1+1)-M_0}{M_0},$

 \hspace{2cm}Res[$\eta^{hom}(1,Y)$, $\pm Bi$] = $-\frac{6(n_2+1)-M_0}{M_0},$

 \hspace{2cm}Res[$\eta^{hom}(1,Y)$, $A \pm Bi$] = $-\frac{6(n_3+1)-M_0}{M_0},$

    where $M_0$ = $2(n_1+1)+ 2(n_2+1)+ 2(n_3+1)$.

These equations are inconsistent so, system (31) is
integrable only when $\mu(x, y) = 0$ i.e.
$\mu_0=\mu_1=\mu_2=\mu_3=\mu_4=\mu_5=0$.

From Proposition 9, system (31) is
axis-reversible if $A=0$ and $\mu_0=\mu_2=\mu_4=0$. Otherwise is
non-reversible because there is no change of variables $\Phi(x, y) = (x+ay, bx+y)$ such that $\Phi_{\ast}F_4$
is axis-reversible, from Theorem 8 the system is non-reversible.
\end{enumerate}
\end{proof}

\vspace{0.5cm}
\textbf{Acknowledgments}\\

This work has been partially supported by Ministerio de Ciencia y Tecnología, Plan Nacional I + D + I co-financed with
FEDER funds, in the frame of the project MTM2010-20907-C02-or, the Consejería de Educación y Ciencia de la Junta de
Andalucía (FQM-276 and PO8-FQM-03770).

\vspace{0.5cm}
\textbf{References}\\

[1] H. Dulac, Détermination et integration d’une certain classe d’équations différentielle ayant par pint singulier un centre, Bull. Sci. Math. Sér. 32 (2)
(1908) 230–252.

[2] N.N. Bautin, On the number of the limit cycles which appear with the variation of coefficients from an equilibrium position of focus or center type,
Mat. Sb. 30 (1952) 181–196; Amer. Math. Soc. Transl. 100 (1954) 1–19.

[3] W. Kapteyn, New investigation on the midpoints of integral of differential equations of the first degree, Nederl. Akad. Wetensch. Verslag Afd. Natuurk.
Konokl. 20 (1912) 1354–1365; Nederl. Akad. Wetensch. Verslag Afd. Natuurk. Konokl. 21 (1912) 27–33 (in Dutch).

[4] W. Kapteyn, On the midpoints of integral curves of differential equations of the first degree, Nederl. Akad. Wetensch. Verslag Afd. Natuurk. Konokl.
Nederland (1911) 1446–1457 (in Dutch).

[5] H. Zoladek, Quadratic system with center and their perturbations, J. Differential Equations 109 (1994) 223–273.

[6] K.S. Sibirskii, On the number of limit cycles on the neighborhood of a singular point, Differ. Equ. 1 (1965) 36–47.

[7] H. Zoladek, On certain generalization of Bautin’s theorem, Nonlinearity (1994) 268–278.

[8] J. Giné, On some of problems in planar differential systems and Hilbert’s 16th problem, Chaos Solitons Fractals 31 (2007) 1118–1134.

[9] M. Berthier, R. Moussu, Reversibilité ét classification des centres nilpotents, Ann. Inst. Fourier (Grenoble) 44 (1994) 465–494.

[10] R. Moussu, Symétrie et forme normale des centres et foyers dégénerés., Ergodic Theory Dynam. Systems 2 (1982) 241–251.

[11] H. Giacomini, J. Giné, J. Llibre, The problem of distinguish between a center or a focus for nilpotents and degenerate analytical system, J. Differential
Equations 227 (2) (2006) 406–426. Corrigendum: J. Differential Equations, 232 (2) (2007) 707.

[12] A. Gassull, J. Torregrosa, Center problem for several differential equations via Cherka’s method, J. Math. Anal. Appl. 228 (1998) 322–342.

[13] A.P. Sadovskii, Problem of distinguishing a center and a focus for a system with a nonvanishing linear part, Differ. Uravn. 12 (7) (1976) 1238–1246.
(Translated).

[14] A. Algaba, C. García, M. Reyes, The center problem for a family of system of differential equations having a nilpotent singular point, J. Math. Anal. Appl.
340 (2008) 32–43.

[15] A. Gassull, J. Llibre, V. Mañosa, F. Mañosas, The focus-center problem for a type of degenerate system, Nonlinearity 13 (2000) 699–730.

[16] A. Algaba, C. García, M.A. Teixeira, Reversibility and quasi-homogeneous normal form of vector fields, Nonlinear Anal. 73 (2) (2010) 510–525.

[17] M.A. Teixeira, J. Yang, The center-focus problem and reversibility, J. Differential Equations 174 (2001) 237–251.

[18] J. Chavarriga, I. García, J. Giné, Integrability of centers perturbed by quasi-homogeneous polynomials, J. Math. Anal. Appl. 210 (1997) 722–730.

[19] J. Chavarriga, H. Giacomini, J. Giné, J. Llibre, Local analytic integrability for nilpotent centers, Ergodic Theory Dynam. Systems 23 (2003) 417–428.

[20] A. Algaba, E. Gamero, C. García, The integrability problem for a class of planar systems, Nonlinearity 22 (2009) 395–420.

[21] L. Cairó, J. Llibre, Polynomial first integrals for weight-homogeneous planar polynomial differential systems of weight degree 3, J. Math. Anal. Appl.
331 (2007) 1284–1298.

[22] J. Llibre, C. Pessoa, On the centers of the weight-homogeneous polynomial vector fields on the plane, J. Math. Anal. Appl. 359 (2009) 722–730.

[23] J. Llibre, X. Zhang, Polynomial first integrals for quasi-homogeneous polynomial differential system, Nonlinearity 15 (2002) 1269–1280.

[24] H. Poincaré, Mémoire sur les courbes dífinies par les équations differentielles, J. Math. Pures Appl. 1 (4) (1885) 167–244.

[25] H. Poincaré, Mémoire sur les courbes dífinies par les équations differentielles, in: Oeuvres de Henri Poincaré, vol. I, Gauthier Villars, París, 1951,
pp. 95–114.

[26] A. Algaba, C. García, M. Reyes, Integrability of two dimensional quasi-homogeneous polynomial differential systems, Rocky Mountain J. Math. 41 (1)
(2011) 1–22.

[27] A. Algaba, E. Freire, E. Gamero, C. García, Monodromy, center-focus and integrability problem for quasi-homogeneous polynomial systems, Nonlinear
Anal. 72 (2010) 1726–1736.

[28] D. Montgomery, L. Zippin, Topological Transformation Groups, Interscience Publ., New York, 1955.
\end{document}